\documentclass{amsart}
\usepackage{amsmath}
\usepackage{amssymb}
\usepackage{graphicx}
\input xy
\xyoption{all}
\newtheorem{Lemma}{Lemma}[section]
\newtheorem{Th}[Lemma]{Theorem}
\newtheorem{Prop}[Lemma]{Proposition}

\newenvironment{Proof}{{\sc Proof.}\ }{~\rule{1ex}{1ex}\vspace{0.2truecm}}
\newcommand{\Cal}[1]{{\mathcal #1}}
\newcommand{\sign}{\mbox{\rm sign}}

\newcommand{\Z}{\mathbb{Z}}

\begin{document}
    \title[On normal subgroups of $D^*$]{On normal subgroups of $D^*$ whose elements are periodic modulo the center of $D^*$ of \\bounded order}
    \dedicatory{Dedicated to Professor Hendrik W. Lenstra for his 65th birthday}
    \author[Mai Hoang Bien]{Mai Hoang Bien}
\address{Mathematisch Instituut, Leiden Universiteit, The Netherlands and Dipartimento di Matematica, Universit\`{a} degli Studi di Padova, Italy.}
\curraddr{Mathematisch Instituut, Leiden Universiteit, Niels Bohrweg 1, 2333 CA Leiden, The Netherlands.}
    \email{maihoangbien012@yahoo.com}

\keywords{Divsion ring, normal subgroup, radical, central. \\
\protect \indent 2010 {\it Mathematics Subject Classification.} 16K20.}

   \maketitle
 \begin{abstract} Let $D$ be a division ring with the center $F=Z(D)$. Suppose that $N$ is a normal subgroup of $D^*$ which is radical over $F$, that is, for any element $x\in N$, there exists a positive integer $n_x$, such that $x^{n_x}\in F$. In \cite{Her1}, Herstein conjectured that $N$ is contained in $F$. In this paper, we show that the conjecture is true  if there exists a positive integer $d$ such that $n_x\le d$ for any $x\in N$.  \end{abstract}

\section{Introduction}
Let $D$ be a division ring with the center $F=Z(D)$. For an element $x\in D$, if there exists a positive integer $n_x$ such that $x^{n_x}\in F$ and $x^{m}\notin F$ for any positive integer $m<n_x$ then $x$ is called {\it $n_x$-central}. If $n_x=1$, $x$ is said to be {\it central}. A subgroup $N$ of the unit group $D^*$ of $D$ is called {\it radical} over $F$ if for any element $x\in N$, there exists $n_x>0$ such that $x$ is $n_x$-central. Such a subgroup $N$ is called {\it central} if $n_x=1$ for any $x\in N$. In other words, $N$ is central if and only if $N$ is contained in $F$. 

 In 1978, Herstein  \cite{Her1} conjectured that if a subnormal subgroup $N$ of $D^*$ is radical over $F$ then it is central. Two years later, he considered the conjecture again and proved that the assumption ``subnormal" in this conjecture is equivalent to ``normal" (see \cite[Lemma 1]{Her2}). That is, he asked whether a normal subgroup of $D^*$  is central  if it is radical over $F$. In \cite{Her1}, Herstein  proved that the conjecture holds if $N$ is torsion. As a consequence, one can see that the conjecture is also true if $D$ is centrally finite. We notice that in \cite{HH}, there is a different proof of this fact. Recall that a division ring $D$ with the center $F$ is called {\it centrally finite} if $D$ is a finite dimensional vector space over $F$ \cite[Definition 14.1]{Lam}. In \cite{Her2}, by using the Pigeon-Hole Principle, Herstein also showed that the conjecture holds if $F$ is uncountable.

Recently, there are some efforts to give the answer for this conjecture. In \cite{HDB1} and \cite{HDB2}, we proved that the conjecture holds if $D$ is either of type $2$ or weakly locally finite. Actually, we get a more general result: if a normal subgroup of $D^*$ is radical over a proper division subring $K$ of $D$ then it is central provided  $D$ is either of type $2$ or weakly locally finite. Recall that a  division ring $D$ is of {\it type $2$} if $\dim_FF(x,y)<\infty$  for any $x, y\in D^*$.  If $F(S)$ is a centrally finite division ring for any finite subset $S$ of  $D$  then $D$ is called {\it weakly locally finite}. Here, $F(S)$ denotes the  division subring of $D$  generated by $F\cup S$. In general, the conjecture remains still open. 

In this paper, we  give a positive answer for  this conjecture in a particular case. In fact, we prove the following Theorem.
\begin{Th} Let $D$ be a division ring and $N$ be a normal subgroup of $D^*$. If there exists a positive integer $d$ such that every element $x\in N$ is $n_x$-central for some positive integer $n_x\le d$ then $N$ is central. \end{Th} 

\section{The proof of the Theorem}
The technique we use in this paper is generalized rational expressions. For our further need, we  recall some definitions and prove some Lemmas.

First, basing on the structure of twisted Laurent series rings, we will construct a division ring which will be used for next Lemmas. Let $R$ be a ring and $\phi$ be a ring automorphism of $R$. We write $\Cal R=R((t,\phi))$ for the ring of formal Laurent series $\sum\limits_{i = n}^\infty  {{a_i}{t^i}}$, where $n\in \Z, a_i\in R$, with the mutiplication defined by the twist equation $ta=\phi(a)t$ for every $a\in R$. In case $\phi(a)=a$ for any $a\in R$, we write $R((t))=R((t,\phi))$. If $R=D$ is a division ring then $\Cal D=D((t,\phi))$ is also a division ring (see \cite[Example 1.8]{Lam}). Moreover, we have.

\begin{Lemma}\label{2.1} Let $R=D$ be a division ring, $\Cal D=D((t,\phi))$ be as above, $F=Z(D)$ be the center of $D$, and $L=\{\, a\in D\mid \phi(a)=a\}$ be the fixed division ring of $\phi$ in $D$. If the center $k=Z(L)$ of $L$ is contained in $F$, then the center of $\Cal D$ is 

$$Z(\Cal D)=\left\{ {\begin{array}{*{20}{c}}
k&{\text{ if } \phi \text{ has infinite order, }}\\
{k(({t^s}))}&{\text{ if } \phi \text{ has an order } s.}
\end{array}} \right.$$ 
\end{Lemma} 
\begin{Proof}  The proof is similar to \cite[Proposition 14.2]{Lam}. It suffices to prove that $Z(D)\subseteq k$ if $\phi$ has infinite order, and $Z(\Cal D)\subseteq k((t^s))$ in case $f$ has an order $s$ since it is easy to check that $k((t^s))\subseteq Z(\Cal D)$ if $f$ has an order $s$. Let $\alpha=\sum\limits_{i = n}^\infty  {{a_i}{t^i}}$ be in $Z(\Cal D)$. We first prove that $a_i\in k$ for every $i\ge n$. One has $\sum\limits_{i = n}^\infty  {{a_i}{t^{i+1}}}=(\sum\limits_{i = n}^\infty  {{a_i}{t^i}})t=t\sum\limits_{i = n}^\infty  {{a_i}{t^i}}=\sum\limits_{i = n}^\infty  {{\phi(a_i)}{t^{i+1}}}$. Hence, $\phi(a_i)=a_i$ for every $i\ge n$. It means $a_i\in L$ for every $i\ge n$. Moreover, for any $a\in L$, $\sum\limits_{i = n}^\infty  {{aa_i}{t^{i}}}=(\sum\limits_{i = n}^\infty  {{a_i}{t^i}})a=\sum\limits_{i = n}^\infty  {{a_i\phi(a)}{t^i}}=\sum\limits_{i = n}^\infty  {{a_ia}{t^{i}}}$. Therefore, $aa_i=a_ia$ for every $i\ge n$. It implies, $a_i\in k$ for every $i\ge n$. Now for any $b\in D$, $\sum\limits_{i = n}^\infty  {{ba_i}{t^{i}}}=(\sum\limits_{i = n}^\infty  {{a_i}{t^i}})b=\sum\limits_{i = n}^\infty  {{a_i\phi^i(b)}{t^i}}=\sum\limits_{i = n}^\infty  {{\phi^i(b)a_i}{t^i}}$, so that $ba_i=\phi^i(b)a_i$ for every $i\ge n$.

{\bf Case 1.} The automorphism $\phi$ has infinite order. For some $i\ne 0$, from the fact that $(b-\phi^i(b))a_i=0$, one has $a_i=0$, which implies $\alpha=a_0\in k$.

{\bf Case 2.} The automorphism $\Phi$ has an order $s$. For any $i$ which is not divided by $n$, since $(b-\phi^i(b))a_i=0$, so that $a_i=0$. Therefore, $\alpha=\sum\limits_{i = m}^\infty  {{a_{si}}{t^{si}}}\in k((t^s))$.
\end{Proof}

Let $\{\,t_i\mid i\in \Z\,\}$ be a countable set of indeterminates and $D$ be a division ring. We construct a family of division rings by the following way. Set $$D_0=D((t_0)), D_1 =D_0((t_1)),$$  $$D_{-1}=D_1((t_{-1})), D_2=D_{-1}((t_{2})),$$ 
   for any $n>1,$ $$ D_{-n}=D_n((t_{-n})),D_{n+1}=D_{-n}((t_{n+1})).$$ Now put $D_{\infty}=\bigcup\limits_{n=-\infty}^{+\infty} {{D_n}}$.  Then $D_\infty$ is a division ring. Assume that $F$ is the center of $D$. By Lemma~\ref{2.1}, it is elementary to prove by induction on $n\ge 0$ that the center of $D_0$ is $F_0=F((t_0))$, the center of $D_{n+1}$ is $F_{n+1}=F_{-n}((t_{n+1}))$ and the center of $D_{-n}$ is $F_{-(n+1)}=F_{n+1}((t_{-(n+1)}))$.  In particular, $F$  is contained in $Z(D_\infty)$. Consider an automorphism $f$ on $D_\infty$ defined by $f(a)=a$ for any $a$ in $D$ and $f(t_i)=t_{i+1}$ for every $i\in \Z$.

\begin{Prop}\label{2.2} Let $D, D_\infty$ and $f$ be as above. Then $\Cal D=D_\infty((t,f))$ is a division ring whose center coincides with the center $F$ of $D$.
\end{Prop}
\begin{Proof} We have $D$ is the fixed division ring of $f$ in $D_\infty$. Since the center $F$ of $D$ is contained in the center of $D_\infty$, $f$ has infinite order and by Lemma~\ref{2.1}, $Z(\Cal D)=F.$\end{Proof}

\bigskip
Recall that a {\it generalized rational expression} of a division ring $D$ is an expression constructed from $D$ and a set of noncommutative indeteminates  using addition, subtraction, multiplication and division. A generalized rational expression over $D$ is called a {\it generalized rational identity} if it vanishes on all permissible substitutions from $D$.  A generalized rational expression $f$ of $D$ is called nontrivial if there exists an extension division ring $D_1$ of $D$ such that $f$ is not a generalized rational identity of $D_1$. The details of  generalized rational identities can be found in \cite{Rowen}.  

Given a positive integer $n$ and $n+1$ noncommutative indeteminates $x,y_1,\cdots, y_n$, put $$g_n(x,y_1,y_2,\cdots, y_n)=\sum\limits_{\delta  \in {S_{n + 1}}} {\sign(\delta ).{x^{\delta (0)}}{y_1}{x^{\delta (1)}}{y_2}{x^{\delta (2)}} \ldots {y_n}{x^{\delta (n)}}}, $$ where $S_{n+1}$ is the symmetric group of $\{\,0,1,\cdots, n\,\}$ and $\sign(\delta)$ is the sign of permutation $\delta$. This is the generalized rational expression defined in \cite{BMM} to connect an algebraic element of degree $n$ and a polynomial. We have the first property of this generalized rational expression.
\begin{Lemma}\label{3.1} Let $D$ be a division ring with the center $F$. For any element $a\in D$, the following are equivalent:
\begin{enumerate}
\item The element $a$ is algebraic over $F$ of degree less than $n$.
\item $g_n(a,r_1,r_2,\cdots, r_n)=0$ for any $r_1, r_2,\cdots, r_n\in D$.
\end{enumerate} 
\end{Lemma}
\begin{Proof} See \cite[Corollary 2.3.8]{BMM}
\end{Proof}

 Let $D$ be a division ring with center $F$ and $a$ be an element of $D$.  Then, by definition, $g_n(axa^{-1}x^{-1},y_1,y_2,\cdots, y_n)$ is also a generalized rational expressions of $D$. Notice that, in general, the expression $g_n(x,y_1,\cdots, y_n)$ is a polynomial but $g_n(axa^{-1}x^{-1},y_1,y_2,\cdots, y_n)$ is not necessary a polynomial. If $a$ is algebraic of degree less than $n$ over $F$ then $g_n(a,y_1,y_2,\cdots, y_n)$ is a trivial generalized rational expression according to Lemma~\ref{3.1}. However, the following Lemma shows that $g_n(axa^{-1}x^{-1},y_1,y_2,\cdots, y_n)$ is always nontrivial if $a$ is not in $F$.

\begin{Lemma}\label{1.1} Let $D$ be a division ring with center $F$. If $a\in D\backslash F$  then the generalized rational expression $g_n(axa^{-1}x^{-1},y_1,y_2,\cdots, y_n)$ is nontrivial. 
\end{Lemma}  
\begin{Proof} Let $D_\infty$, $\Cal D=D_\infty ((t, f))$ and $F$ be as in Proposition~\ref{2.2}. Since $a\notin F$, there exists $c\in D$ such that $c=aba^{-1}b^{-1}\ne 1$. Because $a,b,c$ commute with $t$,  $$(c-1)(1+b^{-1}t)^{-1}+1=a(b+t)a^{-1}(b+t)^{-1}.$$ If $a(b+t)a^{-1}(b+t)^{-1}$ is algebraic over $F$ then so is $(c-1)(1+b^{-1}t)^{-1}$. Hence, $(c-1)^{-1}+b^{-1}(c-1)^{-1}t=((c-1)(1+b^{-1}t)^{-1})^{-1}$ is algebraic over $F$. Let $p(x)=x^m+a_{m-1}x^{m-1}+\cdots +a_1x+a_0$, with $m>0$, be the minimal polynomial of $(c-1)^{-1}+b^{-1}(c-1)^{-1}t$ over $F$. It means $$ 0=((c-1)^{-1}+b^{-1}(c-1)^{-1}t)^m+\cdots +a_1((c-1)^{-1}+b^{-1}(c-1)^{-1}t)+a_0.$$ For instance, $(b^{-1}(c-1)^{-1})^m=0$, a contradiction!  Therefore, $a(b+t)a^{-1}(b+t)^{-1}$ is not algebraic over $F$. Using Lemma~\ref{3.1}, we have $$g_n(a(b+t)a^{-1}(b+t)^{-1},r_1,r_2,\cdots, r_n)\ne 0,$$ for some $r_1,r_2,\cdots, r_n\in \Cal D$. This means $g_n(axa^{-1}x^{-1},y_1,y_2,\cdots, y_n)$ is nontrivial.
\end{Proof}

A polynomial identity ring is a ring $R$ with a non-zero polynomial $P$ vanishing on all permissible substitutions from $R$. In this case, $P$ is called {\it polynomial identity} of $R$ or we say that $R$ {\it satisfies} $P$.  There is a well-known result: a division ring is a polynomial identity division ring if and only if it is centrally finite (see \cite[Theorem 6.3.1]{Her3}). We have a similar property for generalized rational identity division rings.
\begin{Lemma}\label{3.3} Let $D$ be a division ring with the center $F$. If there exists a nontrivial generalized rational identity of $D$ then either $D$ is centrally finite or $F$ is finite.
\end{Lemma}
\begin{Proof} See \cite[Theorem 8.2.15]{Rowen}.
\end{Proof}

Now we are ready to prove our Theorem.\\

\noindent

{\bf Proof of Theorem 1.1}\\

Suppose that  $N$ is not contained in $F$. Then, there exists $a\in N\backslash F$. For any $d+1$ elements $r, r_1,r_2,\cdots, r_d$ of $D$ with $r\ne 0$, since $ara^{-1}r^{-1}\in N$ is $n_{a,r}$-central element for some $0<n_{a,r}\le d$, by Lemma~\ref{3.1}, $$g_d(ara^{-1}r^{-1},r_1,r_2,\cdots, r_d)=0.$$ By Lemma~\ref{1.1}, $g_d(axa^{-1}x^{-1},y_1,y_2,\cdots, y_d)$ is a nontrivial generalized rational identity of $D$. Now, in view of  Lemma~\ref{3.3}, either $D$ is centrally finite or $F$ is finite. If $D$ is centrally finite then $N\subseteq F$ by \cite[Theorem 3.1]{HDB1}. If $F$ is finite then $N$ is torsion, so by \cite[Theorem 8]{Her1}, $N\subseteq F$ . Thus, in both cases we have $N\subseteq F$, a  contradiction. 

\subsection*{Acknowledgment} The author is very thankful to the referee for carefully reading the paper and making useful comments.

\end{document}